\documentclass[12pt,leqno]{amsart}
\usepackage{amsmath, amssymb}
\usepackage{graphicx,color,hyperref}
\usepackage{verbatim, enumitem}

\newcommand\cJ{{\mathcal J}}
\newcommand\cL{{\mathcal L}}

\newcommand\cO{{\mathcal O}}

\newcommand\bN{{\mathbb N}}
\newcommand\bC{{\mathbb C}}
\newcommand\bK{{\mathbb K}}
\newcommand\bZ{{\mathbb Z}}
\newcommand\bP{{\mathbb P}}
\newcommand\bQ{{\mathbb Q}}

\def\ii{{\rm i}\kern1pt}

\def\Ker{\mathop{\rm Ker}}

\def\GG{\mathop{\rm GG}\nolimits}
\def\GL{\mathop{\rm GL}\nolimits}

\def\Sing{\mathop{\rm Sing}}

\def\rank{\mathop{\rm rank}}

\def\ECL{\mathop{\rm ECL}\nolimits}
\def\IEL{\mathop{\rm IEL}\nolimits}
\def\pr{\mathop{\rm pr}\nolimits}
\def\Pic{\mathop{\rm Pic}\nolimits}
\def\rank{\mathop{\rm rank}\nolimits}

\def\proof{\noindent{\it Proof.} }

\def\build#1^#2_#3{\mathrel{\mathop{\null#1}\limits^{#2}_{#3}}}
\def\mertorelbar{\vrule width0.6ex height0.65ex depth-0.55ex}
\def\merto{\mathrel{\mertorelbar\kern1.3pt\mertorelbar\kern1.3pt\mertorelbar
    \kern1.3pt\mertorelbar\kern-1ex\raise0.28ex\hbox{${\scriptscriptstyle>}$}}}

\catcode`\@=11
\newdimen\@rrowlength \@rrowlength=6ex
\def\ssrelbar{\vrule width\@rrowlength height0.64ex depth-0.56ex\kern-4pt}
\def\llra#1{\@rrowlength=#1\ssrelbar\rightarrow}
\catcode`\@=12

\def\semidirect{\mathop{\kern2pt\vrule depth-0.3pt height4.3pt 
\kern-2pt\times}\nolimits}

\def\HOMEPAGE{http://www-fourier.ujf-grenoble.fr/\~{}demailly/manuscripts}

\newdimen\plainitemindent \plainitemindent=18pt
\def\plainitem#1{\par\noindent
\hangindent\plainitemindent\hbox to\plainitemindent{#1\hss}\ignorespaces}

\catcode`\@=11
\def\openup{\afterassignment\@penup\dimen@=}
\def\@penup{\advance\lineskip\dimen@
  \advance\baselineskip\dimen@
  \advance\lineskiplimit\dimen@}
\newdimen\jot \jot=3pt
\newskip\plaincentering \plaincentering=0pt plus 1000pt minus 1000pt
\def\ialign{\everycr{}\tabskip\z@skip\halign}
\def\eqalign#1{\null\,\vcenter{\openup\jot\m@th
  \ialign{\strut\hfil$\displaystyle{##}$&$\displaystyle{{}##}$\hfil
      \crcr#1\crcr}}\,}
\newif\ifdt@p
\def\displ@y{\global\dt@ptrue\openup\jot\m@th
  \everycr{\noalign{\ifdt@p \global\dt@pfalse \ifdim\prevdepth>-1000\p@
      \vskip-\lineskiplimit \vskip\normallineskiplimit \fi
      \else \penalty\interdisplaylinepenalty \fi}}}
\def\@lign{\tabskip\z@skip\everycr{}} 
\def\displaylines#1{\displ@y \tabskip\z@skip
  \halign{\hbox to\displaywidth{$\@lign\hfil\displaystyle##\hfil$}\crcr
    #1\crcr}}
\def\eqalignno#1{\displ@y \tabskip\plaincentering
  \halign to\displaywidth{\hfil$\@lign\displaystyle{##}$\tabskip\z@skip
    &$\@lign\displaystyle{{}##}$\hfil\tabskip\plaincentering
    &\llap{$\@lign##$}\tabskip\z@skip\crcr
    #1\crcr}}
\def\leqalignno#1{\displ@y \tabskip\plaincentering
  \halign to\displaywidth{\hfil$\@lign\displaystyle{##}$\tabskip\z@skip
    &$\@lign\displaystyle{{}##}$\hfil\tabskip\plaincentering
    &\kern-\displaywidth\rlap{$\@lign##$}\tabskip\displaywidth\crcr
    #1\crcr}}
\def\plaincases#1{\left\{\,\vcenter{\normalbaselines\m@th
    \ialign{$##\hfil$&\quad##\hfil\crcr#1\crcr}}\right.}
\def\plainmatrix#1{\null\,\vcenter{\normalbaselines\m@th
    \ialign{\hfil$##$\hfil&&\quad\hfil$##$\hfil\crcr
      \mathstrut\crcr\noalign{\kern-\baselineskip}
      #1\crcr\mathstrut\crcr\noalign{\kern-\baselineskip}}}\,}

\catcode`\@=12

\def\dlraw{\mathrel{\rlap{$\longrightarrow$}\kern-1pt\longrightarrow}}
\def\vlra{\mathrel{\smash-}\joinrel\mathrel{\smash-}\joinrel%
\kern-2pt\longrightarrow}
\def\srelbar{\vrule width0.6ex height0.65ex depth-0.55ex}
\def\merto{\mathrel{\srelbar\kern1.3pt\srelbar\kern1.3pt\srelbar
    \kern1.3pt\srelbar\kern-1ex\raise0.28ex\hbox{${\scriptscriptstyle>}$}}}

\newdimen\claimskip \claimskip=7pt
\long\def\claim#1|#2\endclaim
{\removelastskip\vskip\claimskip\noindent{\bf#1.}
{\it\ignorespaces#2}\vskip\claimskip\noindent}

\font\ninerm=cmr9
\font\ninei=cmmi9
\font\ninesy=cmsy9
\font\ninebf=cmbx9
\font\ninett=cmtt9
\font\nineit=cmti9
\font\ninesl=cmsl9

\font\eightrm=cmr8
\font\eighti=cmmi8
\font\eightsy=cmsy8
\font\eightbf=cmbx8

\font\eightit=cmti8
\font\eightsl=cmsl8
\font\sixrm=cmr6
\font\sixi=cmmi6
\font\sixsy=cmsy6
\font\sixbf=cmbx6
\font\fiverm=cmr5
\font\fivei=cmmi5
\font\fivesy=cmsy5
\font\fivebf=cmbx5

\newfam\itfam
\newfam\slfam
\newfam\bffam
\def\eightpoint{\def\rm{\fam0\eightrm}%
\textfont0=\eightrm \scriptfont0=\sixrm \scriptscriptfont0=\fiverm
 \textfont1=\eighti \scriptfont1=\sixi \scriptscriptfont1=\fivei
 \textfont2=\eightsy \scriptfont2=\sixsy \scriptscriptfont2=\fivesy
 \def\it{\fam\itfam\eightit}%
 \textfont\itfam=\eightit
 \def\sl{\fam\slfam\eightsl}%
 \textfont\slfam=\eightsl
 \def\bf{\fam\bffam\eightbf}%
 \textfont\bffam=\eightbf \scriptfont\bffam=\sixbf
 \scriptscriptfont\bffam=\fivebf
 \normalbaselineskip=9pt
 \setbox\strutbox=\hbox{\vrule height7pt depth2pt width0pt}%
 \normalbaselines\rm}

\def\ninepoint{\def\rm{\fam0\ninerm}%
\textfont0=\ninerm \scriptfont0=\sixrm \scriptscriptfont0=\fiverm
 \textfont1=\ninei \scriptfont1=\sixi \scriptscriptfont1=\fivei
 \textfont2=\ninesy \scriptfont2=\sixsy \scriptscriptfont2=\fivesy
 \def\it{\fam\itfam\nineit}%
 \textfont\itfam=\nineit
 \def\sl{\fam\slfam\ninesl}%
 \textfont\slfam=\ninesl
 \def\bf{\fam\bffam\ninebf}%
 \textfont\bffam=\ninebf \scriptfont\bffam=\sixbf
 \scriptscriptfont\bffam=\fivebf
 \normalbaselineskip=11pt
 \setbox\strutbox=\hbox{\vrule height7pt depth2pt width0pt}%
 \normalbaselines\rm}

\def\plainsection#1{\par\vskip .5cm\penalty -100 
\vbox{\noindent{\sc #1}
\vskip 5pt}
\penalty 500}

\def\Bibitem#1&#2&#3&#4&%
{\hangindent=1.66cm\hangafter=1
\noindent\rlap{\hbox{\rm #1}}\kern1.66cm{\rm #2}{\it #3}{\rm #4.}
\vskip5pt}

\def\square{{\hfill \hbox{
\vrule height 1.453ex  width 0.093ex  depth 0ex
\vrule height 1.5ex  width 1.3ex  depth -1.407ex\kern-0.1ex
\vrule height 1.453ex  width 0.093ex  depth 0ex\kern-1.35ex
\vrule height 0.093ex  width 1.3ex  depth 0ex}}}
\def\bigsquare{{\kern-0.3ex\hbox{
\vrule height 1.7ex  width 0.093ex  depth 0ex\kern-0.093ex
\vrule height 1.8ex  width 1.7ex  depth -1.707ex\kern-0.093ex
\vrule height 1.7ex  width 0.093ex  depth 0ex\kern-1.65ex
\vrule height 0.093ex  width 1.6ex  depth 0ex}\kern0.3ex}}
\def\qed{\phantom{~}$\square$\medskip}
\def\smallskip{\vskip 3pt}
\def\medskip{\vskip 5pt}

\title[Towards the Green-Griffiths-Lang conjecture]{Towards the Green-Griffiths-Lang conjecture}

\author{Jean-Pierre Demailly}
\date{January 23, 2015, revised on March 22, 2015}
\begin{document}

\begin{abstract}
The Green-Griffiths-Lang conjecture stipulates that for every projective 
variety $X$ of general type over $\bC$, there exists a proper algebraic 
subvariety of $X$ containing all non constant entire curves $f:\bC\to X$. 
Using the formalism of directed varieties, we prove here that this 
assertion holds true in case $X$ satisfies a strong general type condition
that is related to a certain jet-semistability property of the tangent
bundle~$T_X$. We then give a sufficient criterion for the Kobayashi 
hyperbolicity of an arbitrary directed variety $(X,V)$.
\end{abstract}

\maketitle
\vskip10pt
\hbox to \textwidth{\hfill\it dedicated to the memory of Salah Baouendi}
\vskip20pt

\plainsection{0. Introduction} 

The goal of this paper is to study the Green-Griffiths-Lang conjecture,
as stated in [GG79] and [Lan86]. It is useful to work in a more general 
context and consider the category of directed projective
manifolds (or varieties). Since the basic problems we deal with are 
birationally invariant, the
varieties under consideration can always be replaced by nonsingular
models. A directed projective manifold is a pair $(X,V)$ where $X$ is
a projective manifold equipped with an analytic linear subspace
$V\subset T_X$, i.e.\ a closed irreducible complex analytic subset $V$
of the total space of~$T_X$, such that each fiber $V_x=V\cap T_{X,x}$
is a complex vector space [If $X$ is not irreducible, $V$ should
rather be assumed to be irreducible merely over each component of $X$,
but we will hereafter assume that our varieties are irreducible]. A
morphism $\Phi:(X,V)\to(Y,W)$ in the category of directed manifolds is
an analytic map $\Phi:X\to Y$ such that $\Phi_*V\subset W$. We refer
to the case $V=T_X$ as being the {\it absolute case}, and to the case
$V=T_{X/S}=\Ker d\pi$ for a fibration $\pi:X\to S$, as being the {\it
  relative case}; $V$ may also be taken to be the tangent space to the
leaves of a singular analytic foliation on~$X$, or maybe even a non
integrable linear subspace of $T_X$.

We are especially interested in {\it entire curves} that are tangent to $V$,
namely non constant holomorphic morphisms $f:(\bC,T_\bC)\to (X,V)$ of directed manifolds. In the absolute case, these are just arbitrary entire curves $f:\bC\to X$. The Green-Griffiths-Lang conjecture, in its strong form, stipulates 

\claim 0.1. GGL conjecture|Let $X$ be a projective variety of general type. Then there exists a proper algebraic variety $Y\subsetneq X$ such that every entire curve $f:\bC\to X$ satisfies $f(\bC)\subset Y$.
\endclaim

\noindent [The weaker form would state that entire curves are
algebraically degenerate, so that\break
$f(\bC)\subset Y_f\subsetneq X$ where $Y_f$ might depend on $f\,$]. 
The smallest admissible algebraic set $Y\subset X$ is by definition the 
{\it entire curve locus} of $X$, defined as the Zariski closure
$$
\ECL(X)=\overline{\bigcup_f f(\bC)}^{\rm Zar}.\leqno(0.2)
$$
If $X\subset\bP^N_\bC$ is defined over a number field $\bK_0$ (i.e.\ by
polynomial equations with equations with coefficients in $\bK_0$) and
$Y=\ECL(X)$, it is expected that for every number field $\bK\supset\bK_0$
the set of $\bK$-points in $X(\bK)\smallsetminus Y$ is
finite, and that this property characterizes $\ECL(X)$ as the smallest
algebraic subset $Y$ of $X$ that has the above property for 
all $\bK$ ([Lan86]). This conjectural arithmetical statement would
be a vast generalization of the Mordell-Faltings theorem, and is 
one of the strong motivations to study the geometric GGL conjecture
as a first step.

\claim 0.3. Problem (generalized GGL conjecture)|Let $(X,V)$ be a 
projective directed manifold. Find
geometric conditions on $V$ ensuring that all entire curves $f:(\bC,
T_\bC)\to(X,V)$ are contained in a proper algebraic subvariety
$Y\subsetneq X$. Does this hold when $(X,V)$ is of general type,
in the sense that the canonical sheaf $K_V$ is big~$?$
\endclaim

As above, we define the entire curve locus set of a pair $(X,V)$ to be 
the smallest admissible algebraic set $Y\subset X$ in the above problem, i.e.
$$
\ECL(X,V)=\overline{\mathop{\kern-40pt\bigcup}_{f:(\bC,T_\bC)\to(X,V)}f(\bC)}^{\rm Zar}.\leqno(0.4)
$$
We say that $(X,V)$ is {\it Brody hyperbolic} if $\ECL(X,V)=\emptyset\,$; as is 
well-known, this is equivalent to Kobayashi hyperbolicity whenever
$X$ is compact.

In case $V$ has no singularities, the {\it canonical sheaf} $K_V$ is 
defined to be $(\det \cO(V))^*$ where $\cO(V)$ is the sheaf of holomorphic
sections of $V$, but in general this naive definition would not work.
Take for instance a generic pencil of
elliptic curves $\lambda P(z)+\mu Q(z)=0$ of degree $3$ in $\bP_\bC^2$,
and the linear space $V$ consisting of the tangents to the fibers
of the rational map $\bP_\bC^2\merto\bP_\bC^1$
defined by $z\mapsto Q(z)/P(z)$. Then $V$ is given by
$$
0\longrightarrow \cO(V)\longrightarrow\cO(T_{\bP_\bC^2})
\build\llra{8ex}^{PdQ-QdP}_{}\cO_{\bP^2_\bC}(6)\otimes\cJ_S
\longrightarrow 0
$$
where $S=\Sing(V)$ consists of the 9 points
$\{P(z)=0\}\cap\{Q(z)=0\}$, and $\cJ_S$ is the corresponding ideal
sheaf of~$S$. Since $\det\cO(T_{\bP^2})=\cO(3)$, we see that
$(\det(\cO(V))^*=\cO(3)$ is ample, thus Problem 0.3 would not have a
positive answer (all leaves are elliptic or singular rational curves
and thus covered by entire curves). An even more ``degenerate''
example is obtained with a generic pencil of conics, in which case
$(\det(\cO(V))^*=\cO(1)$ and $\# S=4$.

If we want to get a positive answer to Problem 0.3, it is
therefore indispensable to give a definition of $K_V$ that incorporates
in a suitable way the singularities of $V\,;$ this will be done in
Def.~1.1 (see also Prop.~1.2). The goal is then to give a positive answer 
to Problem~0.3 under some possibly more restrictive conditions for the 
pair $(X,V)$. These conditions will be expressed in terms of the tower
of Semple jet bundles
$$
(X_k,V_k)\to(X_{k-1},V_{k-1})\to\ldots\to(X_1,V_1)\to(X_0,V_0):=(X,V)
\leqno(0.5)
$$
which we define more precisely in Section~1, following [Dem95]. 
It is constructed inductively by setting $X_k=P(V_{k-1})$
(projective bundle of {\it lines} of $V_{k-1}$), and all $V_k$ have the 
same rank $r=\rank V$, so that $\dim X_k=n+k(r-1)$ where $n=\dim X$.
Entire curve loci have their counterparts for all stages of the Semple tower,
namely, one can define
$$
\ECL_k(X,V)=\overline{\mathop{\kern-40pt\bigcup}_{f:(\bC,T_\bC)\to(X,V)}f_{[k]}(\bC)}^{\rm Zar}.\leqno(0.6)
$$
where $f_{[k]}:(\bC,T_\bC)\to(X_k,V_k)$ is the $k$-jet of $f$.
These are by definition algebraic subvarieties of $X_k$,
and if we denote by $\pi_{k,\ell}:X_k\to X_\ell$ the natural
projection from $X_k$ to $X_\ell$, $0\le\ell\le k$, we get immediately
$$
\pi_{k,\ell}(\ECL_k(X,V))=\ECL_\ell(X,V),\qquad \ECL_0(X,V)=\ECL(X,V).
\leqno(0.7)
$$
Let $\cO_{X_k}(1)$ be the tautological line bundle over $X_k$ associated
with the projective structure. We define the $k$-stage Green-Griffiths 
locus of $(X,V)$ to be
$$
\GG_k(X,V)=\overline{(X_k\smallsetminus\Delta_k)\cap
\bigcap_{m\in\bN}\left(\hbox{base locus of }\cO_{X_k}(m)\otimes 
\pi_{k,0}^*A^{-1}\right)}
\leqno(0.8)
$$
where $A$ is any ample line bundle on $X$ and $\Delta_k=\bigcup_{2\le \ell\le k}
\pi_{k,\ell}^{-1}(D_\ell$) is the union of ``vertical divisors'' 
(see section~1; 
the vertical divisors play no role and have to be removed in this context). 
Clearly, $\GG_k(X,V)$ does not depend on the choice of~$A$.
The basic vanishing theorem for entire 
curves (cf.\ [GG79], [SY96] and [Dem95]) asserts that every entire 
curve $f:(\bC,T_\bC)\to(X,V)$ satisfies all differential equations
$P(f)=0$ arising from sections 
$P\in H^0(X_k,\cO_{X_k}(m)\otimes\pi_{k,0}^*A^{-1})$, hence
$$
\ECL_k(X,V)\subset \GG_k(X,V).
\leqno(0.9)
$$
(For this, one uses the fact that $f_{[k]}(\bC)$ is not contained in any
component of $\Delta_k$, cf.~[Dem95]). It is therefore natural to define
the global Green-Griffiths locus of $(X,V)$ to be
$$
\GG(X,V)=\smash{\bigcap_{k\in\bN}}\pi_{k,0}\left(\GG_k(X,V)\right).
\leqno(0.10)
$$
By (0.7) and (0.9) we infer that
$$
\ECL(X,V)\subset\GG(X,V).
\leqno(0.11)
$$
The main result of [Dem11] (Theorem~2.37 and Cor.~3.4) implies the 
following useful information:

\claim 0.12. Theorem|Assume that $(X,V)$ is of ``general type'', i.e.\ that 
the canonical sheaf $K_V$ is big on $X$.
Then there exists an integer $k_0$ such that $\GG_k(X,V)$ is a proper algebraic
subset of $X_k$ for $k\ge k_0$ $[\,$though $\pi_{k,0}(\GG_k(X,V))$ might still 
be equal to $X$ for all $k\,]$.
\endclaim

In fact, if $F$ is an invertible sheaf on $X$ such that
$K_V\otimes F$ is big, the probabilistic estimates of [Dem11, Cor.~2.38 and 
Cor.~3.4] produce sections of
$$
\cO_{X_k}(m)\otimes\pi_{k,0}^*\cO\Big(\frac{m}{kr}\Big(
1+\frac{1}{2}+\ldots+\frac{1}{k}\Big)F\Big)\leqno(0.13)
$$
for $m\gg k\gg 1$. The (long and involved) proof uses a curvature computation
and singular holomorphic Morse inequalities to show that the line bundles 
involved in (0.11) are big on $X_k$ for $k\gg 1$. One applies this to 
$F=A^{-1}$ with $A$ ample on $X$ to produce sections and conclude that 
$\GG_k(X,V)\subsetneq X_k$. 

Thanks to (0.11), the GGL conjecture is satisfied whenever 
$\GG(X,V)\subsetneq X$. By [DMR10], this happens for instance in the
absolute case when $X$ is a generic hypersurface of 
degree~\hbox{$d\ge\smash{2^{n^5}}$ in $\bP^{n+1}$} (see also [Pau08]
for better bounds in low dimensions, and [Siu02, Siu04]).
However, as already mentioned in [Lan86], very simple examples show that 
one can have $\GG(X,V)=X$ even when $(X,V)$ is of general type, and 
this already occurs in the absolute case as soon as $\dim X\ge 2$. 
A typical example is a product of directed manifolds
$$(X,V)=(X',V')\times(X'',V''),\qquad
V=\pr^{\prime\,*}V'\oplus\pr^{\prime\prime\,*}V''.\leqno(0.14)$$
The absolute case $V=T_X$, $V'=T_{X'}$, $V''=T_{X''}$ on a product of curves
is the simplest instance. It is then easy to check that $\GG(X,V)=X$, cf.\
(3.2). Diverio and Rousseau [DR13] have given many more such examples, 
including the case of indecomposable varieties $(X,T_X)$, e.g.\ Hilbert modular 
surfaces, or more generally compact quotients of bounded symmetric domains
of rank${}\ge 2$. The problem here is the failure of some sort of stability 
condition that is introduced in Section 3. This leads to a somewhat
technical concept of more manageable directed pairs $(X,V)$ that we call 
{\it strongly of general type}, see Def.~3.1. Our main result can be stated

\claim 0.15. Theorem (partial solution to the generalized GGL conjecture)|Let 
$(X,V)$ be a directed pair that is strongly of general type.
Then the Green-Griffiths-Lang conjecture holds true for~$(X,V)$, namely
$\ECL(X,V)$ is a proper algebraic subvariety of $X$.
\endclaim

The proof proceeds through a complicated induction on $n=\dim X$ and 
$k=\rank V$, which is the main reason why we have to introduce
directed varieties, even in the absolute case. An interesting feature of 
this result is that the conclusion
on $\ECL(X,V)$ is reached without having to know anything about 
the Green-Griffiths locus $\GG(X,V)$, even a posteriori. Nevetherless,
this is not yet enough to confirm the GGL conjecture. Our hope is
that pairs $(X,V)$ that are of general type without being strongly
of general type -- and thus exhibit some sort of ``jet-instability'' -- can
be investigated by different methods, e.g.\ by the diophantine
approximation techniques of McQuillan [McQ98]. However, Theorem~0.15
provides a sufficient criterion for Kobayashi hyperbolicity [Kob70, Kob78],
thanks to the following concept of algebraic 
jet-hyperbolicity.

\claim 0.16. Definition|A directed variety $(X,V)$ will be said to be 
algebraically jet-hyperbolic if the induced directed variety structure
$(Z,W)$ on every irreducible algebraic variety $Z$ of~$X$ such that
$\rank W\ge 1$ has a desingularization that is strongly of general 
type $[$see Sections~$2$ and $4$ for the definition of induced directed 
structures and further details$]$. We also say that a projective
manifold $X$ is algebraically jet-hyperbolic if $(X,T_X)$ is.
\endclaim

In this context, Theorem 0.15 yields the following connection between
algebraic jet-hyperbolicity and the analytic concept of Kobayashi 
hyperbolicity.

\claim 0.17. Theorem|Let $(X,V)$ be a directed variety structure on a
projective manifold $X$. Assume that $(X,V)$ is algebraically jet-hyperbolic.
Then $(X,V)$ is Kobayashi hyperbolic.
\endclaim

I would like to thank Simone Diverio and Erwan Rousseau for very
stimulating discussions on these questions. I am grateful to Mihai
P\u{a}un for an invitation at KIAS (Seoul) in August 2014, during
which further very fruitful exchanges took place, and for his
extremely careful reading of earlier drafts of the manuscript.

\plainsection{1. Semple jet bundles and associated canonical sheaves} 

Let $(X,V)$ be a directed projective manifold and $r=\rank V$, that is, the
dimension of generic fibers.
Then $V$ is actually a holomorphic subbundle of $T_X$ on the complement
$X\smallsetminus \Sing(V)$ of a certain minimal analytic set
$\Sing(V)\subsetneq X$ of codimension${}\ge 2$, called hereafter the 
singular set of~$V$. If $\mu:\smash{\widehat X}\to X$ is a proper modification
(a composition of blow-ups with smooth centers, say), we get a directed
manifold $\smash{(\widehat X,\widehat V)}$ by taking $\smash{\widehat V}$ 
to be the closure of $\mu_*^{-1}(V')$, where $V'=V_{|X'}$ is the
restriction of $V$ over a Zariski open set $X'\subset X\smallsetminus \Sing(V)$ 
such that $\mu:\mu^{-1}(X')\to X'$ is a biholomorphism. We will be 
interested in taking modifications realized by iterated blow-ups of 
certain nonsingular subvarieties of the singular set $\Sing(V)$, so 
as to eventually ``improve'' the singularities of $V\,$; outside of
$\Sing(V)$ the effect of blowing-up will be irrelevant, as one
can see easily. Following [Dem11], the canonical sheaf $K_V$ is defined
as follows.

\claim 1.1. Definition|For any directed pair $(X,V)$ with $X$ nonsingular,
we define $K_V$ to be the rank~$1$ analytic sheaf such that
$$
K_V(U)=\hbox{sheaf of locally bounded sections 
of}~~\cO_X(\Lambda^r V^{\prime *})(U\cap X')
$$
where $r=\hbox{rank}(V)$, $X'=X\smallsetminus \Sing(V)$, $V'=V_{|X'}$,
and ``bounded'' means bounded with respect to a smooth hermitian metric
$h$ on $T_X$.
\endclaim

For $r=0$, one can set $K_V=\cO_X$, but this case is trivial: clearly
$\ECL(X,V)=\emptyset$.  The above definition of $K_V$ may look like an
analytic one, but it can easily be turned into an equivalent algebraic
definition:

\claim 1.2. Proposition|Consider the natural morphism 
$\cO(\Lambda^rT_X^*)\to \cO(\Lambda^r V^*)$ where $r=\rank V$
$[\cO(\Lambda^r V^*)$ being defined here as the quotient of 
$\cO(\Lambda^rT_X^*)$ by $r$-forms that have zero restrictions 
to $\cO(\Lambda^rV^*)$ on $X\smallsetminus \Sing(V)\,]$. The 
bidual $\cL_V=\cO_X(\Lambda^r V^*)^{**}$ is an invertible sheaf,
and our natural morphism can be written
$$
\cO(\Lambda^rT_X^*)\to \cO(\Lambda^rV^*)=\cL_V\otimes\cJ_V\subset \cL_V
\leqno(1.2.1)
$$
where $\cJ_V$ is a certain ideal sheaf of $\cO_X$ whose zero set is
contained in $\Sing(V)$ and the arrow on the left is surjective by 
definition. Then 
$$
K_V=\cL_V\otimes \overline{\cJ}_V\leqno(1.2.2)
$$
where $\overline{\cJ}_V$ is the integral closure of $\cJ_V$ in $\cO_X$.
In particular, $K_V$ is always a coherent sheaf.
\endclaim

\proof Let $(u_k)$ be a set of generators of $\cO(\Lambda^rV^*)$ obtained
(say) as the images of a basis $(dz_I)_{|I|=r}$ of $\Lambda^rT^*_X$
in some local coordinates near a point $x\in X$. Write $u_k=g_k\ell$ 
where $\ell$ is a local generator of $\cL_V$ at~$x$. Then $\cJ_V=(g_k)$ by
definition. The boundedness condition expressed in Def.~1.1 means that we 
take sections of the form $f\ell$ where $f$ is a holomorphic function 
on~$U\cap X'$ (and $U$ a neighborhood of $x$), such that 
$$
|f|\le C\sum|g_k|\leqno(1.2.3)
$$
for some constant $C>0$. But then $f$ extends holomorphically to $U$ into 
a function that lies in the integral closure $\overline\cJ_V$, and the latter
is actually characterized analytically by condition~(1.2.3). This proves
Prop.~1.2.\qed

By blowing-up $\cJ_V$ and taking a desingularization $\widehat X$, one 
can always find a {\it log-resolution} of $\cJ_V$ (or $K_V$), i.e.\ a 
modification $\mu:\smash{\widehat X}\to X$ such 
that $\mu^*\cJ_V\subset\cO_{\widehat X}$ is an invertible ideal
sheaf (hence integrally closed); it~follows that
$\mu^*\overline{\cJ}_V=\mu^*\cJ_V$ and
$\mu^*K_V=\mu^*\cL_V\otimes\mu^*\cJ_V$ are invertible sheaves
on~$\smash{\widehat X}$. Notice that for any modification 
$\mu':(X',V')\to (X,V)$, there is always a well defined
natural morphism
$$
\mu^{\prime\,*}K_V\to K_{V'}\leqno(1.3)
$$
(though it need not be an isomorphism, and $K_{V'}$ is possibly non
invertible even when $\mu'$ is taken to be a log-resolution of
$K_V$). Indeed $(\mu')_*=d\mu':V'\to \mu^*V$ is continuous
with respect to ambient hermitian metrics on $X$ and $X'$, and going
to the duals reverses the arrows while preserving boundedness with
respect to the metrics. If $\mu'':X''\to X'$ provides a simultaneous
log-resolution of $K_{V'}$ and $\mu^{\prime\,*}K_V$, we get a non trivial
morphism of invertible sheaves
$$
(\mu'\circ\mu'')^*K_V= \mu^{\prime\prime\,*}\mu^{\prime\,*}K_V
\longrightarrow \mu^{\prime\prime\,*}K_{V'},\leqno(1.4)
$$
hence the bigness of $\mu^{\prime\,*}K_V$ with imply that of 
$\mu^{\prime\prime\,*}K_{V'}$. This is a general principle that we would like
to refer to as the ``monotonicity principle'' for canonical sheaves:
one always get more sections by going to a higher level through a
(holomorphic) modification.

\claim 1.5. Definition|We say that the rank $1$ sheaf
$K_V$ is ``big'' if the invertible sheaf $\mu^*K_V$ is big 
in the usual sense for any log resolution 
$\mu:\smash{\widehat X}\to X$ of~$K_V$. 
Finally, we say that $(X,V)$ is of {\it general type} if there exists
a modification $\mu':(X',V')\to (X,V)$
such that $K_{\smash{V'}}$ is big$\;;$ any higher blow-up 
$\mu'':(X'',V'')\to (X',V')$ then also yields a big canonical sheaf
by $(1.3)$.
\endclaim

Clearly, ``general type'' is a birationally (or bimeromorphically)
invariant concept, by the very definition.  When
$\dim X=n$ and $V\subset T_X$ is a subbundle of rank $r\ge 1$,
one constructs a tower of ``Semple $k$-jet bundles''
$\pi_{k,k-1}:(X_k,V_k)\to (X_{k-1},V_{k-1})$ that are 
$\bP^{r-1}$-bundles, with \hbox{$\dim X_k=n+k(r-1)$} and
$\rank(V_k)=r$. For this, we take $(X_0,V_0)=(X,V)$, and for every 
$k\ge 1$, we set inductively $X_k:=P(V_{k-1})$ and
$$
V_k:=(\pi_{k,k-1})_*^{-1}\cO_{X_k}(-1)\subset T_{X_k},
$$
where $\cO_{X_k}(1)$ is the tautological line bundle on $X_k$,
$\pi_{k,k-1}:X_k=P(V_{k-1})\to X_{k-1}$ the natural projection and
$(\pi_{k,k-1})_*=d\pi_{k,k-1}:T_{X_k}\to
\pi^*_{k,k-1}T_{X_{k-1}}$ its differential (cf.\ [Dem95]). In other terms,
we have exact sequences
$$
\leqalignno{
&0\longrightarrow T_{X_k/X_{k-1}}\longrightarrow V_k
\mathop{\longrightarrow}\limits^{(\pi_{k,k-1})_*}\cO_{X_k}(-1)
\longrightarrow 0,&(1.6)\cr
&0\longrightarrow \cO_{X_k}\longrightarrow(\pi_{k,k-1})^*V_{k-1}\otimes
\cO_{X_k}(1)
\longrightarrow T_{X_k/X_{k-1}}\longrightarrow 0,&(1.7)\cr}
$$
where the last line is the Euler exact sequence associated with the 
relative tangent bundle of $P(V_{k-1})\to X_{k-1}$. Notice that we
by definition of the tautological line bundle we have
$$
\cO_{X_k}(-1)\subset\pi_{k,k-1}^*V_{k-1}
\subset\pi_{k,k-1}^*T_{X_{k-1}},
$$
and also $\rank(V_k)=r$. Let us recall also that for $k\ge 2$, there 
are ``vertical divisors''
$D_k=P(T_{X_{k-1}/X_{k-2}})\subset P(V_{k-1})=X_k$, and that
$D_k$ is the zero divisor of the section of $\cO_{X_k}(1)\otimes
\pi_{k,k-1}^*\cO_{X_{k-1}}(-1)$ induced by the second arrow of the
first exact sequence (1.6), when $k$ is replaced by $k-1$. 
This yields in particular
$$
\cO_{X_k}(1)=\pi_{k,k-1}^*\cO_{X_{k-1}}(1)\otimes\cO(D_k).\leqno(1.8)
$$
By composing the projections we get for all pairs of indices 
$0\le j\le k$ natural morphisms
$$
\pi_{k,j}:X_k\to X_j,\quad
(\pi_{k,j})_*=(d\pi_{k,j})_{|V_k}:V_k\to(\pi_{k,j})^*V_j,
$$
and for every $k$-tuple ${\bf a}=(a_1,\ldots,a_k)\in\bZ^k$
we define
$$
\cO_{X_k}({\bf a})=\bigotimes_{1\le j\le k}\pi_{k,j}^*\cO_{X_j}(a_j),\quad
\pi_{k,j}:X_k\to X_j.
$$
We extend this definition to all weights ${\bf a}\in\bQ^k$ to get
a $\bQ$-line bundle in $\Pic(X)\otimes_\bZ\bQ$. Now, Formula (1.8) yields
$$
\cO_{X_k}({\bf a})=\cO_{X_k}(m)\otimes \cO(-{\bf b}\cdot D)\quad
\hbox{where~~$m=|{\bf a}|=\sum a_j$, ${\bf b}=(0,b_2,\ldots,b_k)$}
\leqno(1.9)
$$
and $b_j=a_1+\ldots+a_{j-1}$, $2\le j\le k$.

When $\Sing(V)\neq\emptyset$, one can always
define $X_k$ and $V_k$ to be the respective closures of $X'_k$, $V'_k$ 
associated with $X'=X\smallsetminus \Sing(V)$ and $V'=V_{|X'}$, where
the closure is taken in the nonsingular ``absolute''
Semple tower $(X^a_k,V^a_k)$ obtained from $(X^a_0,V^a_0)=(X,T_X)$. 
We leave the reader check the following easy (but important) observation.

\claim 1.10. Fonctoriality|If $\Phi:(X,V)\to(Y,W)$ is a morphism of
directed varieties such that $\Phi_*:T_X\to\Phi^*T_Y$ is injective
$($i.e.\ $\Phi$ is an immersion$\,)$, then there is a corresponding
natural morphism $\Phi_{[k]}:(X_k,V_k)\to(Y_k,W_k)$
at the level of Semple bundles. If one merely assumes that the differential
$\Phi_*:V\to\Phi^*W$ is non zero, there is still a well defined
meromorphic map $\Phi_{[k]}:(X_k,V_k)\merto(Y_k,W_k)$ for all~$k\ge 0$.
\endclaim

In case $V$ is singular, the $k$-th Semple bundle $X_k$ will also
be singular, but we can still replace $(X_k,V_k)$ by a suitable modification 
$(\smash{\widehat X}_k,\smash{\widehat V}_k)$ if we want to work with 
a nonsingular model $\smash{\widehat X}_k$ of~$X_k$. The exceptional
set of $\smash{\widehat X}_k$ over $X_k$ can be chosen to lie above
$\Sing(V)\subset X$, and proceeding inductively with respect to $k$, we can 
also arrange the modifications in such a way that we get a tower 
structure $\smash{({\widehat X}_{k+1},
{\widehat V}_{k+1})}\to\smash{({\widehat X}_k,{\widehat V}_k)}\,$; however,
in general, it will not be possible to achieve that 
$\smash{{\widehat V}_k}$ is a subbundle of $T_{\widehat X_k}$.

It is not true that
$K_{\smash{\widehat V}_k}$ is big in case $(X,V)$ is of general type (especially
since the fibers of $X_k\to X$ are towers of $\bP^{r-1}$ bundles, and the
canonical bundles of projective spaces are always negative~!). However,
a twisted version holds true, that can be seen as another instance
of the ``monotonicity principle'' when going to higher stages
in the Semple tower.

\claim 1.11. Lemma|If $(X,V)$ is of general type, then there is a
modification $\smash{(\widehat X,\widehat V)}$ such that all pairs
$\smash{(\widehat X_k,\widehat V_k)}$ of the associated Semple tower
have a twisted canonical bundle $\smash{K_{\widehat V_k}\otimes
\cO_{\widehat X_k}(p)}$ that is still big when one multiplies $K_{\widehat V_k}$ by 
a suitable $\bQ$-line bundle $\cO_{\widehat X_k}(p)$, $p\in\bQ_+$.
\endclaim

\proof First assume that $V$ has no singularities.
The exact sequences (1.6) and (1.7) provide
$$
K_{V_k}:=\det V_k^*=\det(T^*_{X_k/X_{k-1}})\otimes\cO_{X_k}(1)=
\pi^*_{k,k-1}K_{V_{k-1}}\otimes\cO_{X_k}(-(r-1))
$$
where $r=\hbox{rank}(V)$. Inductively we get
$$
K_{V_k}=\pi_{k,0}^*K_V\otimes\cO_{X_k}(-(r-1){\bf 1}),\qquad
{\bf 1}=(1,...,1)\in\bN^k.\leqno(1.11.1)
$$
We know by [Dem95] that $\cO_{X_k}({\bf c})$ is relatively ample over
$X$ when we take the special weight
${\bf c}=(2\,3^{k-2},...,2\,3^{k-j-1},...,6,2,1)$, hence
$$
K_{V_k}\otimes\cO_{X_k}((r-1){\bf 1}+\varepsilon{\bf c})
=\pi_{k,0}^* K_V\otimes\cO_{X_k}(\varepsilon{\bf c})
$$
is big over $X_k$ for any sufficiently small positive rational number 
$\varepsilon\in\bQ^*_+$. Thanks to Formula~(1.9), we can in fact
replace the weight $(r-1){\bf 1}+\varepsilon{\bf c}$ by its total degree
$p=(r-1)k+\varepsilon|{\bf c}|\in\bQ_+$. The general case of a singular 
linear space follows by considering suitable ``sufficiently high'' 
modifications $\smash{\widehat X}$ of $X$, the related directed structure
$\smash{\widehat V}$ on $\smash{\widehat X}$, and embedding 
$\smash{(\widehat X_k,\widehat V_k)}$ in the absolute Semple tower
$\smash{(\widehat X^a_k,\widehat V^a_k)}$ of $\smash{\widehat X}$.
We still have a well defined morphism of rank $1$ sheaves
$$
\pi_{k,0}^*K_{\widehat V}\otimes\cO_{\widehat X_k}(-(r-1){\bf 1})\to
K_{\widehat V_k}
\leqno(1.11.2)
$$
because the multiplier ideal sheaves involved at each stage behave according
to the monoto\-nicity principle applied to the projections
$\smash{\pi^a_{k,k-1}:\widehat X^a_k\to \widehat X^a_{k-1}}$ and their
differentials $(\pi^a_{k,k-1})_*$, which yield well-defined transposed
morphisms from the $(k-1)$-st stage to the $k$-th stage at the level of
exterior differential forms. Our contention follows.\qed

\plainsection{2. Induced directed structure on a subvariety
of a jet space}

Let $Z$ be an irreducible algebraic subset of some $k$-jet bundle
$X_k$ over~$X$, $k\ge 0$. We define the linear subspace 
\hbox{$W\subset T_Z\subset T_{X_k|Z}$} to be the closure
$$
W:=\overline{T_{Z'}\cap V_k}\leqno(2.1)
$$
taken on a suitable Zariski open set $Z'\subset Z_{\rm reg}$ where 
the intersection $T_{Z'}\cap V_k$ has constant rank and 
is a subbundle of $T_{Z'}$. Alternatively, we could also take 
$W$ to be the closure of $T_{Z'}\cap V_k$ in the $k$-th stage
$(X^a_k,V^a_k)$ of the absolute Semple tower, which has the advantage 
of being nonsingular. We say that $(Z,W)$ is the 
{\it induced} directed variety structure; this concept of induced
structure already applies of course in the case $k=0$. 
If $f:(\bC,T_{\bC})\to(X,V)$ is such that $f_{[k]}(\bC)\subset Z$,
then 
$$\hbox{either}~~f_{[k]}(\bC)\subset Z_\alpha\quad\hbox{or}\quad
f'_{[k]}(\bC)\subset W,\leqno(2.2)$$
where $Z_\alpha$ is one of the connected components of $Z\smallsetminus Z'$
and $Z'$ is chosen as in (2.1); especially, if $W=0$, we conclude that
$f_{[k]}(\bC)$ must be contained in one of the $Z_\alpha$'s.
In the sequel, we always consider such a subvariety $Z$ of $X_k$ as a 
directed pair $(Z,W)$ by
taking the induced structure described above. By (2.2), if we proceed
by induction on $\dim Z$, the study of curves tangent to $V$ that have
a $k$-lift $f_{[k]}(\bC)\subset Z$ is reduced to the study of
curves tangent to $(Z,W)$. Let us first quote the 
following easy observation.

\claim 2.3. Observation|For $k\ge 1$, let 
$Z\subsetneq X_k$ be an irreducible algebraic subset that projects
onto~$X_{k-1}$, i.e.\ $\pi_{k,k-1}(Z)=X_{k-1}$.
Then the induced directed variety \hbox{$(Z,W)\subset(X_k,V_k)$},
satisfies
$$1\le \rank W<r:=\rank(V_k).$$
\endclaim

\proof Take a Zariski open subset $Z'\subset Z_{\rm reg}$
such that $W'=T_{Z'}\cap V_k$ is a vector bundle over $Z'$. Since 
$X_k\to X_{k-1}$ is a $\bP^{r-1}$-bundle, $Z$ has codimension at most 
$r-1$ in $X_k$. Therefore $\rank W\ge \rank V_k-(r-1)\ge 1$.
On the other hand, if we had $\rank W=\rank V_k$ generically, then
$T_{Z'}$ would contain $V_{k|Z'}$, in particular it would contain all
vertical directions $T_{X_k/X_{k-1}}\subset V_k$ that are tangent
to the fibers of $X_k\to X_{k-1}$. By taking the flow along vertical 
vector fields, we would
conclude that $Z'$ is a union of fibers of $X_k\to X_{k-1}$ up to an
algebraic set of smaller dimension, but this is excluded since $Z$
projects onto $X_{k-1}$ and $Z\subsetneq X_k$.\qed

\claim 2.4. Definition|For $k\ge 1$, let 
$Z\subset X_k$ be an irreducible algebraic subset of $X_k$. We 
assume moreover that 
$Z\not\subset D_k=P(T_{X_{k-1}/X_{k-2}})$ 
$($and~put here $D_1=\emptyset$ in what follows to avoid to have to 
single out the case $k=1)$. In this situation
we say that $(Z,W)$ is of general type modulo $X_k\to X$ if either
$W=0$, or $\rank W\ge 1$ and there exists 
$p\in\bQ_+$ such that $K_W\otimes\cO_{X_k}(p)_{|Z}$ is big over~$Z$, 
possibly after
replacing $Z$ by a suitable nonsingular model $\smash{\widehat Z}$
$($and pulling-back $W$ and $\cO_{X_k}(p)_{|Z}$ to the nonsingular
variety $\smash{\widehat Z}\,)$.
\endclaim

The main result of [Dem11] mentioned in the introduction as 
Theorem~0.12 implies the following important ``induction step''.

\claim 2.5. Proposition|Let $(X,V)$ be a directed pair where
$X$ is projective algebraic. Take an irreducible algebraic subset
$Z\not\subset D_k$ of the associated $k$-jet Semple bundle
$X_k$ that projects onto~$X_{k-1}$, $k\ge 1$, and assume that the 
induced directed space $(Z,W)\subset(X_k,V_k)$ is of general type 
modulo $X_k\to X$, $\rank W\ge 1$. Then there exists a divisor 
$\Sigma\subset Z_\ell$
in a sufficiently high stage of the Semple tower $(Z_\ell,W_\ell)$
associated with $(Z,W)$, such that every non constant holomorphic map
$f:\bC\to X$ tangent to $V$ that satisfies $f_{[k]}(\bC)\subset Z$
also satisfies $f_{[k+\ell]}(\bC)\subset\Sigma$.
\endclaim

\proof Let $E\subset Z$ be a divisor containing 
$Z_{\rm sing}\cup(Z\cap \pi_{k,0}^{-1}(\Sing(V)))$, chosen so that
on the nonsingular Zariski open set $Z'=Z\smallsetminus E$ all linear spaces
$T_{Z'}$, $V_{k|Z'}$ and $W'=T_{Z'}\cap V_k$ are subbundles of $T_{X_k|Z'}$,
the first two having a transverse intersection on $Z'$.
By taking closures over $Z'$ in the
absolute Semple tower of $X$, we get (singular) directed pairs 
$(Z_\ell,W_\ell)\subset (X_{k+\ell},V_{k+\ell})$, which we eventually
resolve into 
$(\smash{\widehat Z}_\ell,\smash{\widehat W}_{\ell})\subset 
(\smash{\widehat X}_{k+\ell},\smash{\widehat V}_{k+\ell})$ over
nonsingular bases. By construction, locally bounded sections of 
$\cO_{\smash{\widehat X}_{k+\ell}}(m)$ restrict to locally bounded 
sections of $\cO_{\smash{\widehat Z}_{\ell}}(m)$ over 
$\smash{\widehat Z}_{\ell}$.

Since Theorem~0.12 and the related estimate (0.13) are universal in the 
category of directed varieties, we can apply them by replacing $X$
with $\smash{\widehat Z}\subset \smash{\widehat X}_k$, the order $k$ by
a new index $\ell$, and $F$ by 
$$
F_k=\mu^*\Big(\big(\cO_{X_k}(p)\otimes\pi_{k,0}^*\cO_X(-\varepsilon A)
\big)_{|Z}\Big)
$$
where $\mu:\widehat Z\to Z$ is the desingularization, $p\in\bQ_+$ is
chosen such that $K_W\otimes\cO_{x_k}(p)_{|Z}$ is big, 
$A$ is an ample bundle on $X$ and $\varepsilon\in\bQ_+^*$is  small enough.
The assumptions
show that $K_{\widehat{W}}\otimes F_k$ is big on $\smash{\widehat Z}$,
therefore, by applying our theorem and taking $m\gg\ell\gg 1$, 
we get in fine a large number of (metric bounded) sections of
$$\eqalign{
\cO_{\widehat Z_\ell}(m)&\otimes{\widehat \pi}_{k+\ell,k}^*
\cO\Big(\frac{m}{\ell r'}\Big(
1+\frac{1}{2}+\ldots+\frac{1}{\ell}\Big)F_k\Big)\cr
&=
\cO_{{\widehat X}_{k+\ell}}(m{\bf a'})\otimes{\widehat\pi}_{k+\ell,0}^*\cO\Big(
-\frac{m\varepsilon}{kr}\Big(
1+\frac{1}{2}+\ldots+\frac{1}{k}\Big)A\Big)_{|{\widehat Z}_\ell}\cr}
$$
where ${\bf a'}\in\bQ_+^{k+\ell}$ is a positive weight 
(of the form $(0,\ldots,\lambda,\ldots,0,1)$ with some
non zero component $\lambda\in\bQ_+$ at index $k$).
These sections descend to metric bounded sections of
$$
\cO_{X_{k+\ell}}((1+\lambda)m)\otimes{\widehat\pi}_{k+\ell,0}^*\cO\Big(
-\frac{m\varepsilon}{kr}\Big(
1+\frac{1}{2}+\ldots+\frac{1}{k}\Big)A\Big)_{|Z_\ell}.
$$
Since $A$ is ample on $X$, we can apply the fundamental vanishing theorem 
(see e.g.\ [Dem97] or [Dem11], Statement~8.15), or
rather an ``embedded'' version for curves satis\-fying $f_{[k]}(\bC)\subset Z$,
proved exactly by the same arguments. The vanishing theorem
implies that the divisor $\Sigma$ of any such section satisfies the
conclusions of Proposition~2.5, possibly modulo exceptional divisors
of $\smash{\widehat Z}\to Z$; to take care of these, it is enough to add
to $\Sigma$ the inverse image of the divisor $E=Z\smallsetminus Z'$ 
initially selected.\qed

\plainsection{3. Strong general type condition for directed manifolds}

Our main result is the following partial solution to the Green-Griffiths-Lang
conjecture, providing a sufficient algebraic condition for the analytic
conclusion to hold true. We first give an ad hoc definition.

\claim 3.1. Definition|Let $(X,V)$ be a directed pair where
$X$ is projective algebraic. We say that that $(X,V)$ is ``strongly of
general type'' if it is of general type and for every irreducible 
algebraic set $Z\subsetneq X_k$, $Z\not\subset D_k$, that projects 
onto~$X$, the induced directed structure \hbox{$(Z,W)\subset(X_k,V_k)$}
is of general type modulo $X_k\to X$.
\endclaim

\claim 3.2. Example|\rm The situation of a product 
$(X,V)=(X',V')\times(X'',V'')$ described in (0.14) shows that $(X,V)$ can
be of general type without being strongly of general type. In fact,
if $(X',V')$ and $(X'',V'')$ are of general type, then 
$K_V=\pr^{\prime\,*}K_{V'}\otimes\pr^{\prime\prime\,*}K_{V''}$ is big, so $(X,V)$ is 
again of general type. However 
$$
Z=P(\pr^{\prime\,*}V')=X'_1\times X''\subset X_1
$$
has a directed structure $W=\pr^{\prime\,*}V'_1$ which does not possess
a big canonical bundle over $Z$, since the restriction of $K_W$ to
any fiber $\{x'\}\times X''$ is trivial. The higher stages $(Z_k,W_k)$
of the Semple tower of $(Z,W)$ are given by $Z_k=X'_{k+1}\times X''$ and
$W_k=\pr^{\prime\,*}V'_{k+1}$, so it is easy to see that $\GG_k(X,V)$ contains
$Z_{k-1}$. Since $Z_k$ projects onto $X$, we have here $\GG(X,V)=X$
(see [DR13] for more sophisticated indecomposable examples).
\endclaim

\claim 3.3. Remark|\rm It follows from Definition~2.4 that 
$(Z,W)\subset(X_k,V_k)$ is automatically of general type modulo $X_k\to X$ 
if $\cO_{X_k}(1)_{|Z}$ is big. Notice further that
$$
\cO_{X_k}(1+\varepsilon)_{|Z}=
\big(\cO_{X_k}(\varepsilon)\otimes\pi_{k,k-1}^*\cO_{X_{k-1}}(1)
\otimes\cO(D_k)\big)_{|Z}
$$
where $\cO(D_k)_{|Z}$ is effective and $\cO_{X_k}(1)$ is relatively
ample with respect to the projection $X_k\to X_{k-1}$. Therefore the
bigness of $\cO_{X_{k-1}}(1)$ on $X_{k-1}$ also implies that every
directed subvariety $(Z,W)\subset(X_k,V_k)$ is of general type modulo
$X_k\to X$. If
$(X,V)$ is of general type, we know by the main result of [Dem11] that
$\cO_{X_k}(1)$
is big for $k\ge
k_0$ large enough, and actually the precise estimates obtained therein
give explicit bounds for such a~$k_0$.
The above observations show that we need to check the condition of
Definition~3.1 only for $Z\subset
X_k$, $k\le k_0$.  Moreover, at least in the case where
$V$,~$Z$, and $W=T_Z\cap V_k$ are nonsingular, we have
$$
K_W\simeq K_Z\otimes\det (T_Z/W)\simeq K_Z\otimes \det(T_{X_k}/V_k)_{|Z}
\simeq K_{Z/X_{k-1}}\otimes\cO_{X_k}(1)_{|Z}.
$$
Thus we see that, in some sense, it is only needed to check the bigness 
of $K_W$ modulo $X_k\to X$ for ``rather special subvarieties'' $Z\subset X_k$ 
over $X_{k-1}$, such that $K_{Z/X_{k-1}}$ is not relatively big over $X_{k-1}$.\qed
\endclaim

\claim 3.4. Hypersurface case|\rm Assume that $Z\ne D_k$ is an irreducible 
hypersurface of $X_k$ that projects onto $X_{k-1}$. To simplify things further,
also assume that $V$ is nonsingular. Since the Semple jet-bundles $X_k$ 
form a tower of $\bP^{r-1}$-bundles, their Picard groups satisfy
$\Pic(X_k)\simeq\Pic(X)\oplus\bZ^k$ and we have 
$\cO_{X_k}(Z)\simeq\cO_{X_k}({\bf a})\otimes\pi_{k,0}^*B$
for some ${\bf a}\in\bZ^k$ and $B\in\Pic(X)$, where $a_k=d>0$ is the
relative degree of the hypersurface over $X_{k-1}$. Let $\sigma\in H^0(X_k,
\cO_{X_k}(Z))$ be the section defining $Z$ in $X_k$.
The induced directed variety $(Z,W)$ has $\rank W=r-1=\rank V-1$ and
formula (1.12) yields $K_{V_k}=\cO_{X_k}(-(r-1){\bf 1})\otimes\pi_{k,0}^*(K_V)$.
We claim that
$$
K_W\supset\big(K_{V_k}\otimes \cO_{X_k}(Z)\big)_{|Z}
\otimes\cJ_S=
\big(\cO_{X_k}({\bf a}-(r-1){\bf 1})\otimes\pi_{k,0}^*(B\otimes K_V)\big)_{|Z}
\otimes\cJ_S\leqno(3.4.1)
$$
where $S\subsetneq Z$ is the set (containing $Z_{\rm sing}$) where $\sigma$
and $d\sigma_{|V_k}$ both vanish, and $\cJ_S$ is the ideal locally 
generated by the coefficients of $d\sigma_{|V_k}$ along $Z=\sigma^{-1}(0)$. In
fact, the intersection $W=T_Z\cap V_k$ is transverse on $Z\smallsetminus S\,$;
then (3.4.1) can be seen by looking at the morphism
$$
V_{k|Z}\build\llra{4ex}^{d\sigma_{|V_k}}_{}\cO_{X_k}(Z)_{|Z},
$$
and observing that the contraction by $K_{V_k}=\Lambda^rV_k^*$ provides a
metric bounded section of the canonical sheaf $K_W$. In order to investigate
the positivity properties of  $K_W$, one has to show that $B$ cannot be too
negative, and in addition to control the singularity set $S$. The second
point is a priori very challenging, but we get useful information for
the first point by observing that $\sigma$ provides a morphism
$\pi_{k,0}^*\cO_X(-B)\to\cO_{X_k}({\bf a})$, hence a nontrivial morphism
$$
\cO_X(-B)\to E_{\bf a}:=(\pi_{k,0})_*\cO_{X_k}({\bf a})
$$
By [Dem95, Section~12], there exists a filtration on $E_{\bf a}$ such that 
the graded pieces are irreducible representations of $\GL(V)$ contained
in $(V^*)^{\otimes \ell}$, $\ell\le|{\bf a}|$. Therefore we get a nontrivial 
morphism
$$
\cO_X(-B)\to (V^*)^{\otimes \ell},\qquad \ell\le|{\bf a}|.\leqno(3.4.2)
$$
If we know about certain (semi-)stability properties of $V$, this can be used
to control the negativity of $B$.\qed
\endclaim

\noindent
We further need the following useful concept that slightly
generalizes entire curve loci.

\claim 3.5. Definition|If $Z$ is an algebraic set contained
in some stage $X_k$ of the Semple tower of~$(X,V)$, we define its
``induced entire curve locus'' $\IEL_{X,V}(Z)\subset Z$ to be the Zariski closure 
of the union $\bigcup f_{[k]}(\bC)$ of all jets of entire curves
$f:(\bC,T_\bC)\to(X,V)$ such that $f_{[k]}(\bC)\subset Z$.
\endclaim

We~have of course $\IEL_{X,V}(\IEL_{X,V}(Z))=\IEL_{X,V}(Z)$ by definition.
It is not hard to check that modulo certain ``vertical divisors'' of $X_k$, 
the $\IEL_{X,V}(Z)$ locus is essentially the same as 
the entire curve locus $\ECL(Z,W)$ of the induced directed variety,
but we will not use this fact here. Notice that if $Z=\bigcup Z_\alpha$ is
a decomposition of $Z$ into irreducible divisors, then
$$\IEL_{X,V}(Z)=\bigcup_\alpha \IEL_{X,V}(Z_\alpha).$$
Since $\IEL_{X,V}(X_k)=\ECL_k(X,V)$,
proving the Green-Griffiths-Lang property amounts to showing that
$\IEL_{X,V}(X)\subsetneq X$ in the stage $k=0$ of the tower. The basic
step of our approach is expressed in the following statement.

\claim 3.6. Proposition|Let $(X,V)$ be a directed variety and 
$p_0\le n=\dim X$, $p_0\ge 1$.
Assume that there is an integer $k_0\ge 0$ such that for every $k\ge k_0$
and every irreducible algebraic set $Z\subsetneq X_k$, $Z\not\subset D_k$, 
such that $\dim \pi_{k,k_0}(Z)\ge p_0$, the induced directed structure
\hbox{$(Z,W)\subset(X_k,V_k)$} is of general type modulo $X_k\to X$.
Then $\dim\ECL_{k_0}(X,V)<p_0$.
\endclaim

\proof We argue here by contradiction, assuming that 
$\dim\ECL_{k_0}(X,V)\ge p_0$. If 
$$p'_0:=\dim\ECL_{k_0}(X,V)>p_0$$ 
and if we can prove the result for $p'_0$, we will already get a contradiction, 
hence we can assume without loss of generality that $\dim\ECL_{k_0}(X,V)=p_0$.
The main argument
consists of producing inductively an increasing sequence of integers
$$k_0<k_1<\ldots<k_j<\ldots$$
and directed varieties $(Z^j,W^j)\subset(X_{k_j},V_{k_j})$ satisfying
the following properties~:
\vskip3pt
{\plainitemindent=14mm
\plainitem{(3.6.1)} $Z^0$ is one of the irreducible 
components of $\ECL_{k_0}(X,V)$ and $\dim Z^0=p_0\;$;
\vskip3pt
\plainitem{(3.6.2)} $Z^j$ is one of the irreducible 
components of $\ECL_{k_j}(X,V)$ and $\pi_{k_j,k_0}(Z^j)=Z^0\;$;
\vskip3pt
\plainitem{(3.6.3)} for all $j\ge 0$, $\IEL_{X,V}(Z^j)=Z^j$ and
$\rank W_j\ge 1\;$;
\vskip3pt
\plainitem{(3.6.4)} for all $j\ge 0$, the directed variety 
$(Z^{j+1},W^{j+1})$ is contained 
in some stage (of order $\ell_j=k_{j+1}-k_j$) of the Semple
tower of $(Z^j,W^j)$, namely 
$$
(Z^{j+1},W^{j+1})\subsetneq (Z^j_{\ell_j},W^j_{\ell_j})\subset
(X_{k_{j+1}},V_{k_{j+1}})
$$
and
$$
W^{j+1}=\overline{T_{Z^{j+1\,\prime}}\cap W^j_{\ell_j}}=
\overline{T_{Z^{j+1\,\prime}}\cap V_{k_j}}
$$
is the induced directed structure; moreover $\pi_{k_{j+1},k_j}(Z^{j+1})=Z^j$.
\vskip3pt
\plainitem{(3.6.5)} for all $j\ge 0$, we have $Z^{j+1}\subsetneq 
Z^j_{\ell_j}$ but $\pi_{k_{j+1},k_{j+1}-1}(Z^{j+1})=Z^j_{\ell_j-1}$.\vskip3pt}

\noindent
For $j=0$, we simply take $Z^0$ to be one of the irreducible components 
$S_\alpha$ of $\ECL_{k_0}(X,V)$ such that $\dim S_\alpha=p_0$, which exists by 
our hypothesis that $\dim\ECL_{k_0}(X,V)=p_0$. Clearly, $\ECL_{k_0}(X,V)$
is the union of the $\IEL_{X,V}(S_\alpha)$ and we have 
$\IEL_{X,V}(S_\alpha)=S_\alpha$ for all those components, thus 
$\IEL_{X,V}(Z^0)=Z^0$ and $\dim Z^0=p_0$. Assume that $(Z^j,W^j)$
has been constructed. The subvariety $Z^j$ cannot
be contained in the vertical divisor $\smash{D_{k_j}}$. In fact
no irreducible algebraic set $Z$ such that $\IEL_{X,V}(Z)=Z$ can be 
contained in a vertical divisor $D_k$, because $\pi_{k,k-2}(D_k)$ 
corresponds to stationary jets in $X_{k-2}\,$; as every non constant 
curve $f$ has non stationary points, its $k$-jet $f_{[k]}$ cannot
be entirely contained in $D_k\,$; also the induced directed structure
$(Z,W)$ must satisfy $\rank W\ge 1$ otherwise $\IEL_{X,V}(Z)\subsetneq Z$.
Condition (3.6.2) implies that
$\dim\pi_{k_j,k_0}(Z^j)\ge p_0$, thus $(Z^j,W^j)$ is of general type modulo
$X_{k_j}\to X$ by the assumptions of the proposition. Thanks to 
Proposition~2.5, we get an algebraic subset $\smash{\Sigma\subsetneq Z^j_\ell}$
in some stage of the Semple tower $(Z^j_\ell)$ of $Z^j$ such that
every entire curve $f:(\bC,T_\bC)\to(X,V)$ satisfying
$f_{[k_j]}(\bC)\subset Z^j$ also satisfies 
$f_{[k_j+\ell]}(\bC)\subset\Sigma$. By definition, this implies
the first inclusion in the sequence
$$
Z^j=\IEL_{X,V}(Z^j)\subset\pi_{k_j+\ell,k_j}(\IEL_{X,V}(\Sigma))\subset
\pi_{k_j+\ell,k_j}(\Sigma)\subset Z^j
$$
(the other ones being obvious), so we have in fact an equality throughout.
Let $(S'_\alpha)$ be the irreducible
components of $\IEL_{X,V}(\Sigma)$. We have $\IEL_{X,V}(S'_\alpha)=S'_\alpha$ and
one of the components $S'_\alpha$ must satisfy
$$\pi_{k_j+\ell,k_j}(S'_\alpha)=Z^j=Z^j_0.$$
We take $\ell_j\in[1,\ell]$ to be 
the smallest order such that $Z^{j+1}:=\pi_{k_j+\ell,k_j+\ell_j}(S'_\alpha)
\subsetneq\smash{Z^j_{\ell_j}}$, and set $k_{j+1}=k_j+\ell_j>k_j$. 
By definition of $\ell_j$, we have $\smash{\pi_{k_{j+1},k_{j+1}-1}(Z^{j+1})=
Z^j_{\ell_j-1}}$, otherwise $\ell_j$ would not be minimal.
Then $\pi_{k_{j+1},k_j}(Z^{j+1})=Z^j$, hence
$\pi_{k_{j+1},k_0}(Z^{j+1})=Z^0$ by induction, and all properties $(3.6.1-3.6.5)$ 
follow easily. Now, by Observation~2.3, we have
$$
\rank W^j<\rank W^{j-1}<\ldots<\rank W^1<\rank W^0=\rank V.
$$
This is a contradiction because we cannot have such an infinite sequence.
Proposition~3.6 is proved.\qed

\noindent
The special case $k_0=0$, $p_0=n$ of Proposition~3.6 yields the following
consequence.

\claim 3.7. Partial solution to the generalized GGL conjecture|Let 
$(X,V)$ be a directed pair that is strongly of general type.
Then the Green-Griffiths-Lang conjecture holds true for $(X,V)$, namely
$\ECL(X,V)\subsetneq X$, in other words
there exists a proper algebraic variety $Y\subsetneq X$ such that
every non constant holomorphic curve $f:\bC\to X$ tangent to $V$
satisfies $f(\bC)\subset Y$.
\endclaim

\claim 3.8. Remark|\rm The proof is not very constructive, but it is 
however theoretically effective. By this we mean that if $(X,V)$ is
strongly of general type and is taken in a bounded family of 
directed varieties, i.e.\ $X$ is embedded
in some projective space $\bP^N$ with some bound $\delta$ on the degree, 
and  $P(V)$ also has bounded degree${}\le \delta'$ when viewed as 
a subvariety of $P(T_{\bP^N})$, then one could theoretically derive bounds 
$d_Y(n,\delta,\delta')$ for the degree of the locus~$Y$. Also, there 
would exist bounds $k_0(n,\delta,\delta')$ for the orders $k$ and 
bounds $d_k(n,\delta,\delta')$ for the degrees of 
subvarieties $Z\subset X_k$ that have to be checked in the 
definition of a pair of strong general type. In fact, [Dem11] produces
more or less explicit bounds for the order $k$ such that Proposition~2.5
holds true. The degree of the divisor $\Sigma$ is given by a section of
a certain twisted line bundle $\cO_{X_k}(m)\otimes\pi^*_{k,0}\cO_X(-A)$
that we know to be big by an application of holomorphic Morse inequalities 
-- and the bounds for the degrees of $(X_k,V_k)$ then provide bounds
for~$m$.\qed
\endclaim

\claim 3.9. Remark|\rm The condition that $(X,V)$ is strongly of 
general type seems to be related to some sort of stability condition.
We are unsure what is the most appropriate definition, but here is one
that makes sense. Fix an ample divisor $A$ on~$X$. For every 
irreducible subvariety $Z\subset X_k$ that projects onto $X_{k-1}$ for
$k\ge 1$, and $Z=X=X_0$ for $k=0$, we define the slope $\mu_A(Z,W)$
of the corresponding directed variety $(Z,W)$ to be
$$
\mu_A(Z,W)=\frac{\inf\lambda}{\rank W},
$$
where $\lambda$ runs over all rational numbers such that there exists
$m\in\bQ_+$ for which 
$$
K_W\otimes\big(\cO_{X_k}(m)\otimes\pi_{k,0}^*\cO(\lambda A)
\big)_{|Z}\quad\hbox{is big on $Z$}
$$
(again, we assume here that $Z\not\subset D_k$ for $k\ge 2$). Notice that
$(X,V)$ is of general type if and only if $\mu_A(X,V)<0$, and that
$\mu_A(Z,W)=-\infty$ if $\cO_{X_k}(1)_{|A}$ is big. Also, the proof
of Lemma~1.11 shows that 
$$
\mu_A(X_k,V_k)\le \mu_A(X_{k-1},V_{k-1})\le\ldots\le
\mu_A(X,V)\quad\hbox{for all~$k$}
$$
(with $\mu_A(X_k,V_k)=-\infty$ for $k\ge k_0\gg 1$ if $(X,V)$ is of
general type). We say that $(X,V)$ is {\it $A$-jet-stable} (resp.\ 
{\it $A$-jet-semi-stable})
if $\mu_A(Z,W)<\mu_A(X,V)$ (resp.\ $\mu_A(Z,W)\le\mu_A(X,V)$) for
all $Z\subsetneq X_k$ as above. It is then clear that if
$(X,V)$ is of general type and $A$-jet-semi-stable, then it is strongly
of general type in the sense of Definition~3.1. It would be useful 
to have a better understanding of this condition of stability 
(or any other one that would have better properties).\qed
\endclaim

\claim 3.10. Example: case of surfaces|\rm Assume that $X$ is a 
minimal complex surface of general type and $V=T_X$ (absolute case). 
Then $K_X$ is nef and big and the Chern classes of $X$ satisfy
$c_1\le 0$ ($-c_1$ is big and nef) and $c_2\ge 0$. The Semple jet-bundles 
$X_k$ form here a tower of $\bP^1$-bundles and $\dim X_k=k+2$.
Since $\det V^*=K_X$ is big, the strong general type assumption of 
3.6 and 3.8 need only be checked for irreducible hypersurfaces
$Z\subset X_k$ distinct from $D_k$ that project onto $X_{k-1}$, of 
relative degree~$m$. The projection $\pi_{k,k-1}:Z\to X_{k-1}$ is a 
ramified $m:1$ cover. Putting $\cO_{X_k}(Z)\simeq\cO_{X_k}({\bf a})\otimes
\pi_{k,0}(B)$, $B\in\Pic(X)$, we can apply (3.4.1) to get an inclusion
$$
K_W\supset
\big(\cO_{X_k}({\bf a}-{\bf1})\otimes\pi_{k,0}^*(B\otimes K_X)\big)_{|Z}
\otimes\cJ_S,\qquad
{\bf a}\in\bZ^k,~~a_k=m.
$$
Let us assume $k=1$ and $S=\emptyset$ to make things even simpler, and let us
perform numerical calculations in the cohomology ring
$$H^\bullet(X_1,\bZ)=H^\bullet(X)[u]/(u^2+c_1u+c_2),\qquad u=c_1(O_{X_1}(1))
$$
(cf.\ [DEG00, Section~2] for similar calculations and more details). We have
$$
Z\equiv mu+b\quad\hbox{where}\quad b=c_1(B)\quad\hbox{and}\quad
K_W\equiv(m-1)u+b-c_1.$$
We are allowed here to add to $K_W$ an arbitrary multiple 
$\cO_{X_1}(p)$, $p\ge 0$, which we rather write $p=mt+1-m$, $t\ge 1-1/m$.
An evaluation of the Euler-Poincar\'e characteristic of $K_W+\cO_{X_1}(p)_{|Z}$ requires computing the intersection number
$$
\eqalign{
\big(K_W+\cO_{X_1}(p)_{|Z}\big)^2\cdot Z&=\big(mt\,u+b-c_1\big)^2(mu+b)\cr
&=m^2t^2\big(m(c_1^2-c_2)-bc_1\big)
+2mt(b-mc_1)(b-c_1)+m(b-c_1)^2,\cr}
$$
taking into account  that $u^3\cdot X_1=c_1^2-c_2$. In case $S\ne\emptyset$,
there is an additional (negative) contribution from the ideal $\cJ_S$
which is $O(t)$ since $S$ is at most a curve.
In any case, for $t\gg 1$, the leading term in the expansion is 
$m^2t^2(m(c_1^2-c_2)-bc_1)$ and the other terms are negligible with respect
to $t^2$, including the one coming from $S$. We know that $T_X$ is 
semistable with 
respect to $c_1(K_X)=-c_1\ge 0$. Multiplication by the section $\sigma$ yields a
morphism $\pi_{1,0}^*\cO_X(-B)\to\cO_{X_1}(m)$, hence by direct image,
a morphism $\cO_X(-B)\to S^mT^*_X$. Evaluating slopes against $K_X$
(a big nef class), the semistability condition implies 
$bc_1\le \frac{m}{2}c_1^2$, and our leading term is bigger that 
$m^3t^2(\frac{1}{2}c_1^2-c_2)$. We get a positive anwer in the well-known 
case where $c_1^2>2c_2$, corresponding to $T_X$ being almost
ample. Analyzing positivity for the full range of values $(k,m,t)$ and
of singular sets $S$ seems an unsurmountable task at this point; in general,
calculations made in [DEG00] and [McQ99] indicate that the Chern class
and semistability conditions become less demanding for higher order
jets (e.g.\ $c_1^2>c_2$ is enough for $Z\subset X_2$, and 
$c_1^2>\frac{9}{13}c_2$ suffices for $Z\subset X_3$). When $\rank V=1$,
major gains come from the use of Ahlfors currents in combination
with McQuillan's tautological inequalities [McQ98]. We therefore hope for a 
substantial strengthening of the above sufficient conditions, and 
a better understanding of the stability issues, possibly in combination
with a use of  Ahlfors currents and tautological inequalities. In the case
of surfaces, an application of Prop.~3.6 for $k_0=1$ and an analysis
of the behaviour of rank $1$ (multi-)foliations on the surface~$X$
(with the crucial use of [McQ98]) was the main argument used in [DEG00]
to prove the hyperbolicity of very general surfaces of degree $d\ge 21$
in $\bP^3$. For these surfaces, one has $c_1^2<c_2$ and $c_1^2/c_2\to 1$ as
$d\to+\infty$. Applying Prop.~3.6 for higher values $k_0\ge 2$ might 
allow to enlarge the range of tractable surfaces, if the behavior of
rank $1$ (multi)-foliations on $X_{k_0-1}$ can be analyzed independently.
\endclaim

\plainsection{4. Algebraic jet-hyperbolicity implies Kobayashi hyperbolicity}

Let $(X,V)$ be a directed variety, where $X$ is an irreducible projective 
variety; the concept still makes sense when $X$ is singular, by embedding
$(X,V)$ in a projective space $(\bP^N,T_{\bP^N})$ and taking the linear space
$V$ to be an irreducible algebraic subset of $T_{\bP^n}$ that is contained 
in $T_X$ at regular points of~$X$.

\claim 4.1. Definition|Let $(X,V)$ be a directed variety. We say that 
$(X,V)$ is algebraically jet-hyperbolic if
for every $k\ge 0$ and every irreducible algebraic subvariety 
$Z\subset X_k$ that is not contained in the 
union $\Delta_k$ of vertical divisors, the induced directed 
structure $(Z,W)$ either satisfies 
\hbox{$W=0$}, or is of general type modulo $X_k\to X$,
i.e.\ has a desin\-gu\-larization $(\smash{\widehat Z},\smash{\widehat W})$, 
$\mu:\smash{\widehat Z}\to Z$, such that some twisted canonical sheaf
\hbox{$K_{\widehat W}\otimes\mu^*(\cO_{X_k}({\bf a})_{|Z})$}, ${\bf a}\in\bN^k$, 
is big.
\endclaim

\noindent
Proposition~3.6 then gives

\claim 4.2. Theorem|Let $(X,V)$ be an irreducible projective directed 
variety that is algebraically jet-hyperbolic in the sense of the above 
definition. Then $(X,V)$ is Brody $($or Kobayashi$\,)$ hyperbolic, i.e.\ 
$\ECL(X,V)=\emptyset$.
\endclaim

\proof Here we apply Proposition 3.6 with $k_0=0$ and $p_0=1$. It is
enough to deal with subvarieties $Z\subset X_k$ such that 
$\dim\pi_{k,0}(Z)\ge 1$, otherwise $W=0$ and can reduce $Z$ to a smaller
subvariety by (2.2). Then we conclude that
$\dim\ECL(X,V)<1$. All entire curves tangent to $V$ have to be constant, 
and we conclude in fact that $\ECL(X,V)=\emptyset$.\qed
\vskip5mm

\plainsection{References}

{\ninepoint
\parskip 1.5pt plus 0.5pt minus 0.5pt

\Bibitem[Dem95]&Demailly, J.-P.:& Algebraic criteria for Kobayashi
hyperbolic projective varieties and jet differentials.& AMS Summer
School on Algebraic Geometry, Santa Cruz 1995, Proc.\ Symposia in
Pure Math., ed.\ by J.~Koll\'ar and R.~Lazarsfeld, Amer.\ Math.\ Soc., 
Providence, RI (1997), 285–-360&

\Bibitem[Dem97]&Demailly, J.-P.:& Vari\'et\'es hyperboliques et
\'equations diff\'erentielles alg\'e\-bri\-ques.& Gaz.\ Math.\ {\bf 73}
(juillet 1997) 3--23, and {\ninett \HOMEPAGE/cartan\_{}augm.pdf}&

\Bibitem[DEG00]&Demailly, J.-P., El Goul, J.:& Hyperbolicity of
generic surfaces of high degree in projective 3-space.& Amer.\ J.\ 
Math.\ {\bf 122} (2000) 515--546&

\Bibitem[Dem11]&Demailly, J.-P.:& 
Holomorphic Morse inequalities and the Green-Griffiths-Lang
conjecture.& November 2010, arXiv: math.AG/1011.3636, dedicated to the
memory of Ec\-kart Viehweg; Pure and Applied Mathematics 
Quarterly {\bf 7} (2011) 1165--1208&

\Bibitem[DMR10]&Diverio, S., Merker, J., Rousseau, E.:& Effective
algebraic degeneracy.& Invent.\ Math.\ {\bf 180} (2010) 161--223&

\Bibitem[DR13]&Diverio, S., and Rousseau, E.:& The exceptional set and 
the Green-Griffiths locus do not always coincide.& arXiv: math.AG/1302.4756 
(v2)&

\Bibitem[GG79]&Green, M., Griffiths, P.:& Two applications of algebraic
geometry to entire holomorphic mappings.& The Chern Symposium 1979,
Proc.\ Internal.\ Sympos.\ Berkeley, CA, 1979, Springer-Verlag, New York
(1980), 41--74&

\Bibitem[Kob70]&Kobayashi, S.& Hyperbolic manifolds and holomorphic
mappings.& Volume 2 of Pure and Applied Mathematics. Marcel Dekker Inc., 
New York, 1970&

\Bibitem[Kob78]&Kobayashi, S.& Hyperbolic complex spaces.& Volume 318
of Grundlehren der Mathematischen Wissenschaften, Springer-Verlag, 
Berlin, 1998&

\Bibitem[Lan86]&Lang, S.:& Hyperbolic and Diophantine analysis.&
Bull.\ Amer.\ Math.\ Soc.\ {\bf 14} (1986) 159--205&

\Bibitem[McQ98]&McQuillan, M.:& Diophantine approximation and foliations.&
Inst.\ Hautes \'Etudes Sci.\ Publ.\ Math.\ {\bf 87} (1998) 121--174&

\Bibitem[McQ99]&McQuillan, M.:& Holomorphic curves on hyperplane sections 
of $3$-folds.& Geom.\ Funct.\ Anal.\ {\bf 9} (1999) 370--392&

\Bibitem[Pau08]&P\u{a}un, M.:& Vector fields on the total space of 
hypersurfaces in the projective space and hyperbolicity.& Math.\ Ann.\
{\bf 340} (2008) 875--892&

\Bibitem[Siu02]&Siu, Y.T.:& Some recent transcendental techniques in
algebraic and complex geometry.& In: Proceedings of the International
Congress of Mathematicians, Vol.~I (Beijing, 2002), Higher Ed.\ Press, Beijing,
2002, 439--448&

\Bibitem[Siu04]&Siu, Y.T.:& Hyperbolicity in complex geometry.& In: The 
legacy of Niels Henrik Abel, Springer, Berlin, 2004, 543--566&

\Bibitem[SY96]&Siu, Y.T., Yeung, S.K.:& Hyperbolicity of the complement of
a generic smooth curve of high degree in the complex projective plane.& Invent.\
Math.\ {\bf 124} (1996), 573--618&

}

\vskip3mm\noindent
Jean-Pierre Demailly\\
Institut Fourier, Universit\'e Grenoble-Alpes\\
BP74, 100 rue des Maths, 38402 Saint-Martin d'H\`eres, France\\
\emph{e-mail}\/: jean-pierre.demailly@ujf-grenoble.fr

\end{document}